\documentclass[10pt,leqno]{article}

\usepackage{amsmath,amssymb,amsthm,mathrsfs,dsfont}

\usepackage[margin=3cm]{geometry} 

\usepackage{titlesec,hyperref}
\usepackage{extarrows}

\usepackage{color}

\usepackage{fancyhdr}
\titleformat{\subsection}{\it}{\thesubsection.\enspace}{1pt}{}

\newtheorem{theo}{Theorem}[section]
\newtheorem{lemm}[theo]{Lemma}
\newtheorem{defi}[theo]{Definition}

\newtheorem{prop}[theo]{Proposition}
\newtheorem{rema}[theo]{Remark}
\numberwithin{equation}{section}

\allowdisplaybreaks 

\begin{document}
\title{Gevrey regularity and analyticity for Camassa-Holm type systems
\hspace{-4mm}
}
\author{Wei $\mbox{Luo}^1$\footnote{E-mail:  luowei23@mail2.sysu.edu.cn} \quad and\quad
 Zhaoyang $\mbox{Yin}^{1,2}$\footnote{E-mail: mcsyzy@mail.sysu.edu.cn}\\
 $^1\mbox{Department}$ of Mathematics,
Sun Yat-sen University,\\ Guangzhou, 510275, China\\
$^2\mbox{Faculty}$ of Information Technology,\\ Macau University of Science and Technology, Macau, China}

\date{}
\maketitle
\hrule

\begin{abstract}
In this paper we mainly investigate the Cauchy problem of some Camassa-Holm type systems. By constructing a new auxiliary function, we present a generalized Ovsyannikov theorem. By using this theorem and the basic properties of Sobolev-Gevrey spaces, we prove the Gevrey regularity and analyticity of these systems. Moreover, we obtain a lower bound of the lifespan and the continuity of the data-to-solution map.\\

\vspace*{5pt}

\vspace*{5pt}
\end{abstract}

\vspace*{10pt}

\tableofcontents
\section{Introduction}
 In this paper we mainly consider the Cauchy problem for some Camassa-Holm type systems which can be rewritten in the following abstract form:
  \begin{align}
\left\{
\begin{array}{ll}
\frac{du}{dt}=F(t,u(t)),\\
u|_{t=0}=u_0.
\end{array}
\right.
\end{align}
In the following, we will prove the well-posedness of (1.1) in Sobolev-Gevrey spaces under some suitable conditions on the function $F$. The most important and famous equation in (1.1) is the Camassa-Holm equation (CH):
\begin{align}\tag{CH}
\left\{
\begin{array}{ll}
m_t+2m_xu+mu_x=0, ~~m=u-u_{xx},,\\[1ex]
m|_{t=0}=m_0.
\end{array}
\right.
\end{align}
or equivalently
  \begin{align}\tag{CH}\label{CH}
\left\{
\begin{array}{ll}
u_t=-u\partial_x-\partial_x(1-\partial_{xx})^{-1}[u^2+\frac{1}{2}(u_x)^2],\\
u|_{t=0}=u_0.
\end{array}
\right.
\end{align}
$~~~~~~$The Camassa-Holm equation was derived as a model for shallow water waves \cite{Camassa, Constantin.Lannes}. It has been investigated extensively because of its great physical significance in the past two decades. The CH equation has a bi-Hamiltonian structure \cite{Constantin-E,Fokas} and is completely integrable \cite{Camassa,Constantin-P}. The solitary wave solutions of the CH equation were considered in \cite{Camassa,Camassa.Hyman}, where the authors showed that the CH equation possesses peakon solutions of the form $Ce^{-|x-Ct|}$. It is worth mentioning that the peakons are solitons and their shape is alike that of the travelling water waves of greatest height, arising as solutions to the free-boundary problem for incompressible Euler equations over a flat bed (these being the governing equations for water waves),
cf. the discussions in \cite{Constantin2,Constantin.Escher4,Constantin.Escher5,Toland}. Constantin and Strauss verified that the peakon solutions of the CH equation are orbitally stable in \cite{Constantin.Strauss}. \\
  $~~~~~~$The local well-posedness for the CH equation was studied in \cite{Constantin.Escher,Constantin.Escher2,Danchin,Guillermo}. Concretely, for initial profile $u_0\in H^s(\mathbb{R})$ with $s>\frac{3}{2}$, it was shown in \cite{Constantin.Escher,Constantin.Escher2,Guillermo} that the CH equation has a unique solution in $C([0,T);H^s(\mathbb{R}))$. Moveover, the local well-posedness for the CH equation in Besov spaces $C([0,T);B^s_{p,r}(\mathbb{R}))$ with $s>\max(\frac{3}{2},1+\frac{1}{p})$ was proved in \cite{Danchin}. The global existence of strong solutions was established in \cite{Constantin,Constantin.Escher,Constantin.Escher2} under some sign conditions and it was shown in \cite{Constantin,Constantin.Escher,Constantin.Escher2,Constantin.Escher3} that the solutions will blow up in finite time when the slope of initial data was bounded by a negative quantity. The global weak solutions for the CH equation were studied in \cite{Constantin.Molinet} and \cite{Xin.Z.P}. The global conservative and dissipative solutions of CH equation were presented in \cite{Bressan.Constantin} and  \cite{Bressan.Constantin2}, respectively. The analyticity for the solutions of CH equation were investigated in \cite{BHP} and \cite{HM1}.\\
$~~~~~~$A natural idea is to extend such study to the multi-component generalized
systems. One of the most popular generalized systems is the following integrable two-component
Camassa-Holm shallow water system (2CH) \cite{Constantin.Ivanov}:
\begin{align}\tag{2CH}
\left\{
\begin{array}{ll}
m_t+um_x+2u_xm+k\rho\rho_x=0, ~~m=u-u_{xx},\\[1ex]
\rho_t+(u\rho)_x=0,                        \\[1ex]
m|_{t=0}=m_0, \rho|_{t=0}=\rho_0,
\end{array}
\right.
\end{align}
or equivalently
  \begin{align}\tag{2CH}\label{2CH}
\left\{
\begin{array}{ll}
u_t=-u\partial_x-\partial_x(1-\partial_{xx})^{-1}[u^2+\frac{1}{2}(u_x)^2+\frac{k}{2}\rho^2],\\[1ex]
\rho_t=-(u\rho)_x,           \\[1ex]
u|_{t=0}=u_0, \rho|_{t=0}=\rho_0,
\end{array}
\right.
\end{align}
where $k=\pm1$. Local well-posedness for (2CH) with the initial
data in Sobolev spaces and in Besov spaces was established in \cite{Constantin.Ivanov}, \cite{Escher.Yin}, and  \cite{GuiGuilong}, respectively. The blow-up phenomena and global existence of strong solutions to (2CH) in
Sobolev spaces were obtained in \cite{Escher.Yin}, \cite{Guan.Yin} and \cite{GuiGuilong}. The existence
of global weak solutions for (2CH) with $k=1$ was investigated in \cite{Guan.weak}.\\
$~~~~~~$Another one is the modified two-component Camassa-Holm
system (M2CH) \cite{Holm.Naraigh}:
\begin{align}\tag{M2CH}
\left\{
\begin{array}{ll}
m_t+um_x+2u_xm+k\rho\overline{\rho}_x=0,~~ m=u-u_{xx}\\[1ex]
\rho_t+(u\rho)_x=0, ~~\rho=(1-\partial^2_x)(\overline{\rho}-\overline{\rho}_0)\\[1ex]
m|_{t=0}=u_0, \rho|_{t=0}=\rho_0,
\end{array}
\right.
\end{align}
or equivalently
  \begin{align}\tag{M2CH}\label{M2CH}
\left\{
\begin{array}{ll}
u_t=-u\partial_x-\partial_x(1-\partial_{xx})^{-1}[u^2+\frac{1}{2}(u_x)^2+\frac{k}{2}\gamma^2-\frac{k}{2}\gamma^2_x],\\[1ex]
\gamma_t=-u\gamma_x-(1-\partial_{xx})^{-1}((u_x\gamma_x)_x+u_x\gamma),           \\[1ex]
u|_{t=0}=u_0, \gamma|_{t=0}=\gamma_0,
\end{array}
\right.
\end{align}
where $k=\pm1$ and $\overline{\rho}_0$ is a constant. Local well-posedness for (M2CH) with the initial
data in Sobolev spaces and in Besov spaces was established in \cite{Guan.Karlsen} and \cite{Kai.Yin} respectively. The blow up phenomena of strong solutions to (M2CH) were presented in \cite{Guan.Karlsen}. The existence
of global weak solutions for (M2CH) with $k=1$ was investigated in \cite{Guan.weak.modified}.  The global conservative and dissipative solutions of (M2CH) equation were studied in \cite{Tan.Yin} and  \cite{Tan.Yin2}, respectively. The analyticity of the solutions for (M2CH) was proved in \cite{Kai.Yin1}.\\
$~~~~~~$ Recently Geng and Xue proposed a new three-component Camassa-Holm system with N-peakon solutions \cite{Geng.Xue}:
  \begin{align}\tag{3CH}
\left\{
\begin{array}{ll}
u_{t}=-va_x+u_xb+\frac{3}{2}ub_x-\frac{3}{2}u(a_xc_x-ac),\\[1ex]
v_t=2vb_x+v_xb,\\[1ex]
w_t=-vc_x+w_xb+\frac{3}{2}wb_x+\frac{3}{2}w(a_xc_x-ac),\\[1ex]
u=a-a_{xx},\\[1ex]
v=\frac{1}{2}(b_{xx}-4b+a_{xx}c_x-c_{xx}a_x+3a_xc-3ac_x), \\[1ex]
w=c-c_{xx},\\[1ex]
u|_{t=0}=u_{0},~~ v|_{t=0}=v_0,~~ w|_{t=0}=w_0.\\[1ex]
\end{array}
\right.
\end{align}
 It is based on the following spectral problem
\begin{align}
\phi_x=U\phi,~~\phi=
\begin{pmatrix}
\phi_1 \\
\phi_2 \\
\phi_3
\end{pmatrix},~~~
U=
\begin{pmatrix}
0 & 1 & 0 \\
1+\lambda v & 0 & u \\
\lambda w   & 0 & 0
\end{pmatrix},
\end{align}
where $u,~v,~w$ are three potentials and $\lambda$ is a constant spectral parameter. It was shown in \cite{Geng.Xue} that the N-peakon solitons of the system (1.1) have the form
\begin{align}
a(t,x)=\sum^N_{i=0} a_i(t)e^{-|x-x_i(t)|},\\
\nonumber b(t,x)=\sum^N_{i=0} b_i(t)e^{-2|x-x_i(t)|},\\
\nonumber c(t,x)=\sum^N_{i=0} c_i(t)e^{-|x-x_i(t)|},
\end{align}
where $a_i,~b_i,~c_i$ and $x_i$ evolve according to a dynamical system. Moreover, the author derived infinitely many conservation laws of the system (1.1). In \cite{Wei.JDE}, the authors proved the local well-posedness and global existence of strong solution to (3CH) under some sign conditions.\\
$~~~~~~$Many researchers have studied the analyticity of solutions to Camassa-Holm type systems, cf.  \cite{BHP}, \cite{HM1} and \cite{Kai.Yin1}. However, to our best acknowledge, the Gevrey regularity of solutions to the Camassa-Holm equation is still an open problem.
Our motivation is to solve this problem.  To begin with, we introduce an abstract Cauchy-Kovalevsky theorem which is very crucial to study the analyticity:
\begin{theo}\cite{BG1,Nr,O1}\label{T}
Let $\{X_\delta\}_{0<\delta<1}$ be a scale of decreasing Banach spaces, namely,  for any $\delta'<\delta$ we have $X_\delta\subset X_{\delta'}$ and $\|\cdot\|_{\delta'}\leq \|\cdot\|_{\delta}$,
and let $T,~R>0$, $\sigma\geq 1$. For given $u_0\in X_1$, assume that the function $F$ satisfies the following
conditions:\\
$(1)$ If for $0 < \delta' < \delta < 1$ the function $t\mapsto u(t)$ is holomorphic in $|t|<T$ and continuous on $|t| < T$ with values in $X_s$ and
$$\sup_{|t|<T}\|u(t)\|_{\delta}<R,$$
then $t \mapsto F(t,u(t))$ is a holomorphic function on $|t|<T$ with values in $X_{\delta'}$.\\
$(2)$ For any $0 < \delta' < \delta < 1$ and any $u,v\in \overline{B(u_0,R)}\subset X_\delta$, there exists a positive constant $L$ depending on $u_0$ and $R$ such that
$$\sup_{|t|<T}\|F(t,u)-F(t,v)\|_{\delta'}\leq\frac{L}{\delta-\delta'}\|u-v\|_{\delta}.$$
$(3)$ For any $0<\delta<1$, there exists a positive constant $M$ depending on $u_0$ and $R$ such that
$$\sup_{|t|<T}\|F(t,u_0)\|_{\delta}\leq\frac{M}{1-\delta}.$$
Then  there exists a $T_0\in(0,T)$ and a unique solution to the Cauchy problem (1.1), which for every $\delta\in(0,1)$ is holomorphic in $|t|<T_0(1-\delta)$ with values in $X_\delta$.
\end{theo}
  Theorem \ref{T} was first proposed by Ovsyannikov in \cite{O1},\cite{O2},\cite{O3}.
  However, the original Ovsyannikov theorem becomes invalid for the Gevrey class. Because this kind of spaces do not satisfy the condition (2) of the Ovsyannikov theorem. More precisely, in Section 2, for the Gevrey class, we see that
  \begin{align}
  \sup_{|t|<T}\|F(t,u)-F(t,v)\|_{\delta'}\leq\frac{L}{(\delta-\delta')^\sigma}\|u-v\|_{\delta},
  \end{align}
  with $\sigma\geq 1$. If $\sigma>1$, the inequality (1.4) is weaker than the condition (2) because it is nonlinear decay. Thus, we need a new framework which is associated with the properties of the Gevrey class. In this paper, we modify the proof of \cite{O1} and establish a new auxiliary function, then obtain a generalised Ovsyannikov theorem. By using this theorem, we obtain both the Gevrey regularity and analyticity of the solutions to Camassa-Holm type systems. Moreover, by taking advantage of the idea in \cite{BHP}, we prove that the continuity of the data-to-solution map.\\
  $~~~~~~$ The paper is organized as follows. In Section 2 we recall some properties about Sobolev-Gevrey spaces. In Section 3, we prove a generalized Ovsyannikov theorem. In Section 4, we prove the analyticity and Gevrey regularity of the solutions to some Camassa-Holm type systems. In Section 5, we show that the data-to-solution map is continuous from the data space to the solution space.
\section{Preliminaries}
Firstly, we introduce the Sobolev-Gevrey spaces and recall some basic properties.
\begin{defi}\cite{Foias}\label{Gevrey}
Let $s$ be a real number and $\sigma, \delta>0$. A function $f \in G^{\delta}_{\sigma,s}(\mathbb{R}^d)$ if and only if $f\in C^\infty(\mathbb{R}^d)$ and satisfies
$$\|f\|_{G^\delta_{\sigma,s}(\mathbb{R}^d)}=\bigg(\int_{\mathbb{R}^d}(1+|\xi|^2)^se^{2\delta|\xi|^{\frac{1}{\sigma}}}|\widehat{f}(\xi)|^2d\xi\bigg)^{\frac{1}{2}}<\infty.$$
\end{defi}
\begin{rema}
Denoting the Fourier multiplier $e^{\delta(-\Delta)^{\frac{1}{2\sigma}}}$ by
$$e^{\delta(-\Delta)^{\frac{1}{2\sigma}}}f= \mathscr{F}^{-1}(e^{\delta|\xi|^{\frac{1}{\sigma}}}\widehat{f}),$$
we deduce that $\|f\|_{G^\delta_{\sigma,s}(\mathbb{R}^d)}=\|e^{\delta(-\Delta)^{\frac{1}{2\sigma}}}f\|_{H^s(\mathbb{R}^d)}$. For $0<\sigma<1$, it is called ultra-analytic function. If $\sigma=1$, it is usual analytic function and $\delta$ is called the radius of analyticity. If $\sigma>1$, it is the Gevrey class function.
\end{rema}
\begin{prop}\label{p3}
Let $0<\delta'<\delta$, $0<\sigma'<\sigma$ and $s'<s$. From Definition \ref{Gevrey}, one can check that $G^{\delta}_{\sigma,s}(\mathbb{R}^d)\hookrightarrow G^{\delta'}_{\sigma,s}(\mathbb{R}^d)$, $G^{\delta}_{\sigma',s}(\mathbb{R}^d)\hookrightarrow G^{\delta}_{\sigma,s}(\mathbb{R}^d)$ and $G^{\delta}_{\sigma,s}(\mathbb{R}^d)\hookrightarrow G^{\delta}_{\sigma,s'}(\mathbb{R}^d)$.
\end{prop}

\begin{prop}\label{D}
Let $s$ be a real number and $\sigma>0$. Assume that $0<\delta'<\delta$. Then we have
$$\|\partial_{x}f\|_{G^{\delta'}_{\sigma,s}(\mathbb{R})}\leq \frac{e^{-\sigma}\sigma^{\sigma}}{(\delta-\delta')^\sigma}\|f\|_{G^{\delta}_{\sigma,s}(\mathbb{R})}.$$
\begin{proof}
Since $\widehat{\partial_{x}f}=i\xi \widehat{f}$, it follows that
\begin{align}
\|\partial_{x}f\|^2_{G^{\delta'}_{\sigma,s}(\mathbb{R})}&=\int_{\mathbb{R}}(1+|\xi|^2)^se^{2\delta'|\xi|^{\frac{1}{\sigma}}}|\xi|^2|\widehat{f}(\xi)|^2d\xi\\
\nonumber&=\frac{1}{(\delta-\delta')^{2\sigma}}\int_{\mathbb{R}}(1+|\xi|^2)^se^{2\delta|\xi|^{\frac{1}{\sigma}}}e^{-2[(\delta-\delta')^\sigma|\xi|]^{\frac{1}{\sigma}}}(\delta-\delta')^{2\sigma}|\xi|^2|\widehat{f}(\xi)|^2d\xi\\
\nonumber&\leq \frac{\|f\|^2_{G^{\delta}_{\sigma,s}(\mathbb{R})}}{(\delta-\delta')^{2\sigma}} \sup_{\xi\in\mathbb{R}} \{e^{-2[(\delta-\delta')^\sigma|\xi|]^{\frac{1}{\sigma}}}(\delta-\delta')^{2\sigma}|\xi|^2\}.
\end{align}
Let $z=[(\delta-\delta')^\sigma|\xi|]^{\frac{1}{\sigma}}\geq 0$ and consider the function $g(z)=e^{-2z}z^{2\sigma}$. By directly calculating, we have $\lim_{z\rightarrow 0}g(z)=0$, $\lim_{z\rightarrow +\infty}g(z)=0$ and $g'(z)=-2e^{-2z}z^{2\sigma}+2\sigma e^{-2z}z^{2\sigma-1}$. By solving $g'(z)=0$, we obtain that $z=\sigma$, which implies that $g(z)\leq g(\sigma)=e^{-2\sigma}\sigma^{2\sigma}$. Then, we deduce from (2.1) that
\begin{align*}
\|\partial_{x}f\|_{G^{\delta'}_{\sigma,s}(\mathbb{R})}\leq \frac{e^{-\sigma}\sigma^{\sigma}\|f\|_{G^{\delta}_{\sigma,s}(\mathbb{R})}}{(\delta-\delta')^{\sigma}}.
\end{align*}
\end{proof}
\end{prop}

\begin{prop}\label{Product}(Product acts on Sobolev-Gevrey spaces with $d=1$)
Let $s>\frac{1}{2}$, $\sigma\geq 1$ and $\delta>0$. Then, $G^\delta_{\sigma,s}(\mathbb{R})$ is an algebra. Moreover, there exists a constant $C_s$ such that
$$\|fg\|_{G^{\delta}_{\sigma,s}(\mathbb{R})}\leq C_s\|f\|_{G^\delta_{\sigma,s}(\mathbb{R})}\|g\|_{G^\delta_{\sigma,s}(\mathbb{R})}.$$
\begin{proof}
Since $\widehat{fg}=\widehat{f}\ast\widehat{g}$, it follows that
\begin{align}
\|fg\|^2_{G^{\delta}_{\sigma,s}(\mathbb{R})}&=\int_{\mathbb{R}}(1+|\xi|^2)^se^{2\delta|\xi|^{\frac{1}{\sigma}}}|\widehat{f}\ast\widehat{g}|^2d\xi\\
\nonumber&=\int_{\mathbb{R}}(1+|\xi|^2)^s\bigg(\int_{\mathbb{R}}e^{\delta|\xi|^{\frac{1}{\sigma}}}\widehat{f}(\eta)\widehat{g}(\xi-\eta)d\eta\bigg)^2d\xi\\
\nonumber&\leq \int_{\mathbb{R}}(1+|\xi|^2)^s\bigg(\int_{\mathbb{R}}e^{\delta|\xi-\eta|^{\frac{1}{\sigma}}}e^{\delta|\eta|^{\frac{1}{\sigma}}}\widehat{f}(\eta)\widehat{g}(\xi-\eta)d\eta\bigg)^2d\xi \quad (\textit{Here we use the fact that $\sigma\geq 1$})\\
\nonumber&=\int_{\mathbb{R}}(1+|\xi|^2)^s |\mathscr{F}(e^{\delta(-\Delta)^{\frac{1}{2\sigma}}}f)\ast \mathscr{F}(e^{\delta(-\Delta)^{\frac{1}{2\sigma}}}g)|^2 d\xi=\|(e^{\delta(-\Delta)^{\frac{1}{2\sigma}}}f)\cdot (e^{\delta(-\Delta)^{\frac{1}{2\sigma}}}g)\|^2_{H^s(\mathbb{R})}\\
\nonumber&\leq C_s\|e^{\delta(-\Delta)^{\frac{1}{2\sigma}}}f\|^2_{H^s(\mathbb{R})}\|e^{\delta(-\Delta)^{\frac{1}{2\sigma}}}g\|^2_{H^s(\mathbb{R})}\quad (\textit{Here we use the fact that $s>\frac{1}{2}$})\\
\nonumber&=C_s\|f\|^2_{G^{\delta}_{\sigma,s}(\mathbb{R})}\|g\|^2_{G^{\delta}_{\sigma,s}(\mathbb{R})}.
\end{align}
\end{proof}
\end{prop}

\begin{prop}\label{Morse}
Let $s>\frac{1}{2}$, $\sigma\geq 1$ and $\delta>0$. There exists a constant $\overline{C}_s$ such that
$$\|fg\|_{G^{\delta}_{\sigma,s-1}(\mathbb{R})}\leq \overline{C}_s\|f\|_{G^\delta_{\sigma,s-1}(\mathbb{R})}\|g\|_{G^\delta_{\sigma,s}(\mathbb{R})}.$$
\begin{proof}
By the similar argument as in Proposition \ref{Product}, we have
\begin{align}
\|fg\|^2_{G^{\delta}_{\sigma,s-1}(\mathbb{R})}\leq \|(e^{\delta(-\Delta)^{\frac{1}{2\sigma}}}f)\cdot (e^{\delta(-\Delta)^{\frac{1}{2\sigma}}}g)\|^2_{H^{s-1}(\mathbb{R})}.
\end{align}
Using the fact that $\|ab\|_{H^{s-1}(\mathbb{R})}\leq \overline{C}_s\|a\|_{H^{s-1}(\mathbb{R})}\|b\|_{H^{s}(\mathbb{R})}$ if $s>\frac{1}{2}$, we get \begin{align}
\|(e^{\delta(-\Delta)^{\frac{1}{2\sigma}}}f)\cdot (e^{\delta(-\Delta)^{\frac{1}{2\sigma}}}g)\|^2_{H^{s-1}(\mathbb{R})}\leq \overline{C}_s\|e^{\delta(-\Delta)^{\frac{1}{2\sigma}}}f\|^2_{H^{s-1}(\mathbb{R})}\|e^{\delta(-\Delta)^{\frac{1}{2\sigma}}}g\|^2_{H^{s}(\mathbb{R})}=\overline{C}_s\|f\|^2_{G^\delta_{\sigma,s-1}(\mathbb{R})}\|g\|^2_{G^\delta_{\sigma,s}(\mathbb{R})}.
\end{align}
\end{proof}
\end{prop}

Through out this paper, we use the notations $P_1\doteq(1-\partial_{xx})^{-1}$, $P_2\doteq(4-\partial_{xx})^{-1}$, $P_3\doteq\partial_x$ and $P_{ij}\doteq P_iP_j$ with $1\leq i,j\leq 3$.
\begin{prop}\label{Mult}
If $s\in\mathbb{R}$, $\sigma,\delta>0$ and $f\in G^\delta_{\sigma,s}(\mathbb{R})$, then
\begin{align}
&\|P_1f\|_{G^\delta_{\sigma,s}(\mathbb{R})}=\|f\|_{G^\delta_{\sigma,s-2}(\mathbb{R})}\leq\|f\|_{G^\delta_{\sigma,s}(\mathbb{R})},\\
&\|P_2f\|_{G^\delta_{\sigma,s}(\mathbb{R})}\leq \frac{1}{4}\|f\|_{G^\delta_{\sigma,s}(\mathbb{R})}, \\
&\|P_{13}f\|_{G^\delta_{\sigma,s}(\mathbb{R})}\leq \|f\|_{G^\delta_{\sigma,s-1}(\mathbb{R})}, \\
&\|P_{13}f\|_{G^\delta_{\sigma,s}(\mathbb{R})}\leq \frac{1}{2}\|f\|_{G^\delta_{\sigma,s}(\mathbb{R})}, \\
&\|P_{23}f\|_{G^\delta_{\sigma,s}(\mathbb{R})}\leq \frac{1}{4}\|f\|_{G^\delta_{\sigma,s}(\mathbb{R})}.
\end{align}
\begin{proof}
Since $\mathscr{F}[P_1f]=\frac{\widehat{f}(\xi)}{1+|\xi|^2}$, $\mathscr{F}[P_2f]=\frac{\widehat{f}(\xi)}{4+|\xi|^2}$, $\mathscr{F}[P_{13}f]=\frac{i\xi\widehat{f}(\xi)}{1+|\xi|^2}$ and $\mathscr{F}[P_{23}f]=\frac{i\xi\widehat{f}(\xi)}{4+|\xi|^2}$, it follows that
\begin{align}
&\|P_1f\|_{G^\delta_{\sigma,s}(\mathbb{R})}=(\int_{\mathbb{R}}(1+|\xi|^2)^{s-2}e^{2\delta|\xi|}|\widehat{f}(\xi)|^2d\xi)^{\frac{1}{2}}= \|f\|_{G^\delta_{\sigma,s-2}(\mathbb{R})}\leq \|f\|_{G^\delta_{\sigma,s}(\mathbb{R})},\\
\nonumber&\|P_2f\|_{G^\delta_{\sigma,s}(\mathbb{R})}=(\int_{\mathbb{R}}\frac{(1+|\xi|^2)^{s}}{(4+|\xi|^2)^2}e^{2\delta|\xi|}|\widehat{f}(\xi)|^2d\xi)^{\frac{1}{2}}\leq \frac{1}{4}(\int_{\mathbb{R}}(1+|\xi|^2)^{s}e^{2\delta|\xi|}|\widehat{f}(\xi)|^2d\xi)^{\frac{1}{2}}=\frac{1}{4}\|f\|_{G^\delta_{\sigma,s}(\mathbb{R})},\\
\nonumber&\|P_{13}f\|_{G^\delta_{\sigma,s}(\mathbb{R})}=(\int_{\mathbb{R}}(1+|\xi|^2)^{s-2}|\xi|^2e^{2\delta|\xi|}|\widehat{f}(\xi)|^2d\xi)^{\frac{1}{2}}\leq (\int_{\mathbb{R}}(1+|\xi|^2)^{s-1}e^{2\delta|\xi|}|\widehat{f}(\xi)|^2d\xi)^{\frac{1}{2}}=\|f\|_{G^\delta_{\sigma,s-1}(\mathbb{R})},\\
\nonumber&\|P_{13}f\|_{G^\delta_{\sigma,s}(\mathbb{R})}=(\int_{\mathbb{R}}(1+|\xi|^2)^{s-2}|\xi|^2e^{2\delta|\xi|}|\widehat{f}(\xi)|^2d\xi)^{\frac{1}{2}}\leq (\int_{\mathbb{R}}\frac{(1+|\xi|^2)^{s}}{4}e^{2\delta|\xi|}|\widehat{f}(\xi)|^2d\xi)^{\frac{1}{2}}=\frac{1}{2}\|f\|_{G^\delta_{\sigma,s}(\mathbb{R})},\\
\nonumber&\|P_{23}f\|_{G^\delta_{\sigma,s}(\mathbb{R})}=(\int_{\mathbb{R}}\frac{(1+|\xi|^2)^{s}|\xi|^2}{(4+|\xi|^2)^2}e^{2\delta|\xi|}|\widehat{f}(\xi)|^2d\xi)^{\frac{1}{2}}\leq \frac{1}{4}(\int_{\mathbb{R}}(1+|\xi|^2)^{s}e^{2\delta|\xi|}|\widehat{f}(\xi)|^2d\xi)^{\frac{1}{2}}=\frac{1}{4}\|f\|_{G^\delta_{\sigma,s}(\mathbb{R})}.
\end{align}
\end{proof}
\end{prop}
{\bf Notations}. Since all function spaces in the following sections are over $\mathbb{R}$, for simplicity, we drop $\mathbb{R}$ in the notation of function spaces if there is no ambiguity.
\section{A generalized Ovsyannikov theorem}
In order to study the Gevrey regularity of (1.1), we need the following generalized Ovsyannikov theorem.
\begin{theo}\label{Th1}
Let $\{X_\delta\}_{0<\delta<1}$ be a scale of decreasing Banach spaces, namely,  for any $\delta'<\delta$ we have $X_\delta\subset X_{\delta'}$ and $\|\cdot\|_{\delta'}\leq \|\cdot\|_{\delta}$. Consider the
Cauchy problem
\begin{align}\label{Cauchy}
\left\{
\begin{array}{ll}
\frac{du}{dt}=F(t,u(t)),\\
u|_{t=0}=u_0.
\end{array}
\right.
\end{align}
Let $T,~R>0$, $\sigma\geq 1$. For given $u_0\in X_1$, assume that $F$ satisfies the following
conditions:\\
$(1)$ If for $0 < \delta' < \delta < 1$ the function $t\mapsto u(t)$ is holomorphic in $|t|<T$ and continuous on $|t| < T$ with values in $X_s$ and
$$\sup_{|t|<T}\|u(t)\|_{\delta}<R,$$
then $t \mapsto F(t,u(t))$ is a holomorphic function on $|t|<T$ with values in $X_{\delta'}$.\\
$(2)$ For any $0 < \delta' < \delta < 1$ and any $u,v\in \overline{B(u_0,R)}\subset X_\delta$, there exists a positive constant $L$ depending on $u_0$ and $R$ such that
$$\sup_{|t|<T}\|F(t,u)-F(t,v)\|_{\delta'}\leq\frac{L}{(\delta-\delta')^\sigma}\|u-v\|_{\delta}.$$
$(3)$ For any $0<\delta<1$, there exists a positive constant $M$ depending on $u_0$ and $R$ such that
$$\sup_{|t|<T}\|F(t,u_0)\|_{\delta}\leq\frac{M}{(1-\delta)^\sigma}.$$
Then  there exists a $T_0\in(0,T)$ and a unique solution $u(t)$ to the Cauchy problem (\ref{Cauchy}), which for every $\delta\in(0,1)$ is holomorphic in $|t|<\frac{T_0(1-\delta)^\sigma}{2^\sigma-1}$ with values in $X_\delta$.
\end{theo}
\begin{rema}
In fact, $T_0=\min\{\frac{1}{2^{2\sigma+4}L}, \frac{(2^\sigma-1)R}{(2^\sigma-1)2^{2\sigma+3}LR+M}\}$, which gives a lower bound of the lifespan.
\end{rema}
\begin{rema}
If $\sigma=1$, Theorem \ref{Th1} reduced to the so called abstract Cauchy-Kovalevsky theorem. The original results were first proposed by Ovsyannikov in \cite{O1},\cite{O2} and \cite{O3}. Later, Nirenberg \cite{Nr}, Nishida \cite{Ns}, Treves \cite{Tre1}, \cite{Tre2}, and Baouendi and Goulaouic \cite{BG1}, \cite{BG2} developed a lot of different versions of this theorem.
\end{rema}
The proof of Theorem \ref{Th1} is based on the fixed point argument in some suitable Banach space. Now we introduce a new Banach space.
\begin{defi}\label{Def}
Let $\sigma\geq1$. For any $a>0$ we denote by $E_a$ the space of functions $u(t)$ which for every $0<\delta<1$ and $|t|<\frac{a(1-\delta)^\sigma}{2^\sigma-1}$, are holomorphic and continuous functions of $t$ with values in $X_\delta$ such that
\begin{align}
\|u\|_{E_a}=\sup_{|t|<\frac{a(1-\delta)^\sigma}{2^\sigma-1},0<\delta<1}\bigg(\|u(t)\|_{\delta}(1-\delta)^\sigma\sqrt{1-\frac{|t|}{a(1-\delta)^\sigma}}\bigg)<+\infty.
\end{align}
\end{defi}
\begin{prop}
Let $\sigma\geq 1$. For any $a>0$, the function space $E_a$ is a Banach space equipped with the norm $\|\cdot\|_{E_a}$.
\begin{proof}
Suppose that $(u_n)_{n\geq1}$ is a Cauchy sequence in $E_a$, that is
$$\|u_n-u_m\|_{E_a}\rightarrow0, \quad \text{as} \quad n,m\rightarrow\infty.$$
By virtue of the definition of $E_a$, we deduce that for any $0<\delta<1$,
$$\sup_{|t|<\frac{a(1-\delta)^\sigma}{2^\sigma-1}}\|u_n-u_m\|_\delta\rightarrow 0, \quad \text{as} \quad n,m\rightarrow\infty.$$
Since $X_\delta$ is a Banach space, it follows that there exists a $u_\delta\in X_\delta$ such that
$$\sup_{|t|<\frac{a(1-\delta)^\sigma}{2^\sigma-1}}\|u_n-u_\delta\|_\delta \rightarrow 0, \quad \text{as} \quad n\rightarrow\infty. $$
Now we claim that $u_\delta$ is independent on $\delta$. Indeed, if $\delta_1\neq\delta_2$, with loss of generality suppose that $\delta_1<\delta_2$, and we obtain that,
$$\|u_n-u_{\delta_2}\|_{\delta_1}\leq \|u_n-u_{\delta_2}\|_{\delta_2}\rightarrow 0, \quad \text{as} \quad n\rightarrow\infty,$$
which leads to $u_{\delta_1}=u_{\delta_2}$. Thus, for any $0<\delta<1$, we have $u=u_\delta\in X_\delta$. Since $(u_n)_{n\geq1}$ is a Cauchy sequence in $E_a$, for any $\varepsilon>0$, there exists a $N_1=N_1(\varepsilon)$ such that if $n,m\geq N_1$, $\|u_n-u_m\|_{E_a}\leq \frac{\varepsilon}{2}$. Note that $\|u_n-u\|_{\delta}\xrightarrow{n\rightarrow\infty}0$ for any $0<\delta<1$. For any $\varepsilon>0$, there exists a $N_2(\delta)$ such that if $n\geq N_2(\delta)$, $\|u_n-u\|_{\delta}\leq \frac{\varepsilon}{2}$. Defining that $N=N(\delta,\varepsilon)=\max\{N_1,N_2(\delta)\}+1$ for any $\varepsilon>0$ and $0<\delta<1$, we deduce that for any $n\geq N_1$
\begin{align*}
\|u_n-u\|_\delta(1-\delta)^\sigma\sqrt{1-\frac{|t|}{a(1-\delta)^\sigma}}&\leq \|u_n-u_N\|_{E_a}+\|u_N-u\|_\delta(1-\delta)^\sigma\sqrt{1-\frac{|t|}{a(1-\delta)^\sigma}}\\
&\leq \|u_n-u_N\|_{E_a}+\|u_N-u\|_\delta\leq \frac{\varepsilon}{2}+\frac{\varepsilon}{2}=\varepsilon.
\end{align*}
Since $N_1$ is independent on $\delta$, it follows from the above inequality  that $\|u_n-u\|_{E_a}\xrightarrow{n\rightarrow\infty}0$.
\end{proof}
\end{prop}
The following lemmas are crucial to prove Theorem \ref{Th1}.
\begin{lemm}\label{L.1}
Let $\sigma\geq1$. For every $0<\delta<1$ and $0\leq t<\frac{a(1-\delta)^\sigma}{2^\sigma-1}$ we have
$$1-\delta> (\frac{1}{2})^{1+\frac{1}{\sigma}}\bigg\{[(1-\delta)^\sigma-\frac{t}{a}]^{\frac{1}{\sigma}}+[(1-\delta)^\sigma+(2^{\sigma+1}-1)\frac{t}{a}]^{\frac{1}{\sigma}}\bigg\}.$$
\begin{proof}
Since $t<\frac{a(1-\delta)^\sigma}{2^\sigma-1}$, it follows that
\begin{align}
2(1-\delta)^\sigma > (1-\delta)^\sigma+(2^\sigma-1)\frac{t}{a}= \frac{1}{2}[(1-\delta)^\sigma-\frac{t}{a}]+\frac{1}{2}[(1-\delta)^\sigma+(2^{\sigma+1}-1)\frac{t}{a}].
\end{align}
Using the fact that $(x+y)^p\leq 2^{p-1}(x^p+y^p)$ with $p=\sigma$, $x=(\frac{1}{2}[(1-\delta)^\sigma-\frac{t}{a}])^\frac{1}{\sigma}>0$ and $y=(\frac{1}{2}[(1-\delta)^\sigma+(2^{\sigma+1}-1)\frac{t}{a}])^\frac{1}{\sigma}>0$, we deduce that
\begin{align}
\frac{1}{2}[(1-\delta)^\sigma-\frac{t}{a}]+\frac{1}{2}[(1-\delta)^\sigma+(2^{\sigma+1}-1)\frac{t}{a}]=x^\sigma+y^\sigma\geq \frac{(x+y)^\sigma}{2^{\sigma-1}}.
\end{align}
Plugging (3.4) into (3.3) yields that
\begin{align}
1-\delta>\frac{x+y}{2}=(\frac{1}{2})^{1+\frac{1}{\sigma}}\bigg\{[(1-\delta)^\sigma-\frac{t}{a}]^{\frac{1}{\sigma}}+[(1-\delta)^\sigma+(2^{\sigma+1}-1)\frac{t}{a}]^{\frac{1}{\sigma}}\bigg\}.
\end{align}
\end{proof}
\end{lemm}
\begin{lemm}\label{L.2}
Let $\sigma\geq 1$. For every $a>0$, $u\in E_a$, $0<\delta<1$ and $0\leq t<\frac{a(1-\delta)^\sigma}{2^\sigma-1}$ we have
$$\int^t_0\frac{\|u(\tau)\|_{\delta(\tau)}}{(\delta(\tau)-\delta)^\sigma}d\tau\leq \frac{a2^{2\sigma+3}\|u\|_{E_a}}{(1-\delta)^\sigma}\sqrt{\frac{a(1-\delta)^\sigma}{a(1-\delta)^\sigma-t}},$$
where $\delta(\tau)=\frac{1}{2}(1+\delta)+(\frac{1}{2})^{2+\frac{1}{\sigma}}\bigg\{[(1-\delta)^\sigma-\frac{t}{a}]^\frac{1}{\sigma}-[(1-\delta)^\sigma+(2^{\sigma+1}-1)\frac{t}{a}]^\frac{1}{\sigma}\bigg\}\in (\delta,1)$.
\begin{proof}
By virtue of the definition of $E_a$, we obtain
\begin{align}
\int^t_0\frac{\|u(\tau)\|_{\delta(\tau)}}{(\delta(\tau)-\delta)^\sigma}d\tau\leq\|u\|_{E_a}\int^t_0\frac{1}{(\delta(\tau)-\delta)^\sigma(1-\delta(\tau))^\sigma\sqrt{1-\frac{\tau}{a(1-\delta(\tau))^\sigma}}}d\tau.
\end{align}
Taking advantage of Lemma \ref{L.1}, we have
\begin{align}
\delta(\tau)-\delta&=\frac{1}{2}(1-\delta)+(\frac{1}{2})^{2+\frac{1}{\sigma}}\bigg\{[(1-\delta)^\sigma-\frac{\tau}{a}]^\frac{1}{\sigma}-[(1-\delta)^\sigma+(2^{\sigma+1}-1)\frac{\tau}{a}]^\frac{1}{\sigma}\bigg\}\\
\nonumber&\geq   (\frac{1}{2})^{1+\frac{1}{\sigma}}[(1-\delta)^\sigma-\frac{\tau}{a}]^\frac{1}{\sigma},
\end{align}
and
\begin{align}
1-\delta(\tau)& =\frac{1}{2}(1-\delta) -(\frac{1}{2})^{2+\frac{1}{\sigma}}\bigg\{[(1-\delta)^\sigma-\frac{\tau}{a}]^\frac{1}{\sigma}-[(1-\delta)^\sigma+(2^{\sigma+1}-1)\frac{\tau}{a}]^\frac{1}{\sigma}\bigg\}
\\ \nonumber&\geq(\frac{1}{2})^{1+\frac{1}{\sigma}} [(1-\delta)^\sigma+(2^{\sigma+1}-1)\frac{\tau}{a}]^\frac{1}{\sigma},
\end{align}
which leads to
\begin{align}
(1-\delta(\tau))^\sigma\geq (\frac{1}{2})^{\sigma+1}[(1-\delta)^\sigma-\frac{\tau}{a}]+\frac{\tau}{a},
\end{align}
or equivalently
\begin{align}
a(1-\delta(\tau))^\sigma-\tau \geq (\frac{1}{2})^{\sigma+1}[a(1-\delta)^\sigma-\tau].
\end{align}
Plugging (3.7)-(3.10) into (3.6) yields that
\begin{align}
\int^t_0\frac{\|u(\tau)\|_{\delta(\tau)}}{(\delta(\tau)-\delta)^\sigma}d\tau&\leq \|u\|_{E_a}\int^t_0\frac{a^2}{[a(1-\delta)^{\sigma}-\tau]^{\frac{3}{2}}[a(1-\delta)^{\sigma}+(2^{\sigma+1}-1)\tau]^{\frac{1}{2}}}d\tau\\
\nonumber&=\frac{a2^{2(\sigma+1)}}{(1-\delta)^\sigma}\|u\|_{E_a}\int^{\frac{t}{a(1-\delta)^{\sigma}}}_0\frac{1}{(1-\theta)^\frac{3}{2}(1+(2^{\sigma+1}-1)\theta)^\frac{1}{2}}d\theta\\
\nonumber&\leq \frac{a2^{2(\sigma+1)}}{(1-\delta)^\sigma} \|u\|_{E_a}\int^{\frac{t}{a(1-\delta)^{\sigma}}}_0\frac{1}{(1-\theta)^\frac{3}{2}} d\theta\leq \frac{a2^{2\sigma+3}\|u\|_{E_a}}{(1-\delta)^\sigma}\sqrt{\frac{a(1-\delta)^\sigma}{a(1-\delta)^\sigma-t}}.
\end{align}
\end{proof}
\end{lemm}

 \textit{Proof of Theorem \ref{Th1}:}
 We only consider the case $t\geq 0$. For any $t<\frac{a(1-\delta)^\sigma}{2^\sigma-1}$ with $a>0$ and $u(t)\in \overline{B(u_0,R)}\subset E_a$, we define that
\begin{align}
G(u(t))\doteq u_0+\int^t_0F(\tau,u(\tau))d\tau.
\end{align}
Since (3.1) is equivalent to
\begin{equation}
u(t)=u_0+\int^t_0F(\tau,u(\tau))d\tau,
\end{equation}
it follows that our initial value problem (3.1) can be reduced to find the fixed point of the operator $G$. \\
{\bf Step 1}: If $u(t)\in E_a$, by virtue of Definition \ref{Def}, we have $u(t)$ is a holomorphic and continuous function of $t$ with values in $X_\delta$ for any $0<\delta<1$. The condition (1) of $F$ implies that $F(t,u(t))$ is a holomorphic  function of $t$ with values in $X_\delta$ for any $0<\delta<1$, which leads to $G(u(t))$ is a holomorphic and continuous function of $t$ with values in $X_\delta$ for any $0<\delta<1$. In addition, if $\|u-u_0\|_{E_a}\leq R$, we deduce from Lemma \ref{L.2} and conditions (2)-(3) that
\begin{align}
\|G(u(t))-u_0\|_{\delta}&\leq \int^t_0\|F(\tau,u(\tau))\|_{\delta}d\tau\leq \int^t_0\|F(\tau,u(\tau))-F(\tau,u_0)\|_{\delta}d\tau+\int^t_0\|F(\tau,u_0)\|_{\delta}d\tau\\
\nonumber&\leq \int^t_0\frac{L\|u-u_0\|_{\delta(\tau)}}{(\delta(\tau)-\delta)^\sigma}d\tau+ \frac{tM}{(1-\delta)^\sigma}\leq \frac{a2^{2\sigma+3}LR}{(1-\delta)^\sigma}\sqrt{\frac{a(1-\delta)^\sigma}{a(1-\delta)^\sigma-t}}+\frac{tM}{(1-\delta)^\sigma},
\end{align}
which implies that
\begin{align}
\|G(u(t))-u_0\|_{E_a}\leq a2^{2\sigma+3}LR+\frac{aM}{2^\sigma-1}.
\end{align}
By taking $a\leq \frac{(2^\sigma-1)R}{(2^\sigma-1)2^{2\sigma+3}LR+M}$, we verify that $Gu\in \overline{B(u_0,R)}\subset E_a$, which leads to $G$ maps $\overline{B(u_0,R)}\subset E_a$ into itself. \\
{\bf Step 2}: Assume that $u(t),v(t)\in \overline{B(u_0,R)}\subset E_a$. Taking advantage of Lemma \ref{L.2} and the condition (2), we infer that
\begin{align}
\|G(u(t))-G(v(t))\|_{\delta}&\leq \int^t_0\|F(\tau,u(\tau))-F(\tau,v(\tau))\|_{\delta}d\tau\leq \\
\nonumber&\leq \int^t_0\frac{L\|u-v\|_{\delta(\tau)}}{(\delta(\tau)-\delta)^\sigma}d\tau\leq \frac{a2^{2\sigma+3}L\|u-v\|_{E_a}}{(1-\delta)^\sigma}\sqrt{\frac{a(1-\delta)^\sigma}{a(1-\delta)^\sigma-t}},
\end{align}
which leads to
\begin{align}
\|G(u(t))-G(v(t))\|_{E_a}\leq a2^{2\sigma+3}L\|u-v\|_{E_a}.
\end{align}
By taking $a\leq \frac{1}{2^{2\sigma+4}L}$, we obtain $\|G(u(t))-G(v(t))\|_{E_a}\leq \frac{1}{2}\|u-v\|_{E_a}$, and hence $G$ is a contraction map on $\overline{B(u_0,R)}\subset E_a$. From Step 1 and Step 2, we deduce that if $a\leq T_0=\min\{\frac{1}{2^{2\sigma+4}L}, \frac{(2^\sigma-1)R}{(2^\sigma-1)2^{2\sigma+3}LR+M}\}$, $T$ has a unique fixed point in $\overline{B(u_0,R)}\subset E_a$.
\section{Gevrey regularity and analyticity}
In this section we investigate the Gevrey regularity and analyticity of solutions to the Camassa-Holm type systems. By virtue of Remark 2.2, the case $\sigma>1$ is corresponding to Gevrey regularity while $\sigma=1$ is corresponding to analyticity. Our main results can be stated as follows.
\begin{theo}\label{CHTH}
Let $\sigma\geq1$ and $s>\frac{3}{2}$. Assume that $u_0\in G^1_{\sigma,s}(\mathbb{R})$. Then for every $0<\delta<1$,  there exists a $T_0>0$ such that the Camassa-Holm equation has a unique solution $u$ which is holomorphic in $|t|<\frac{T_0(1-\delta)^\sigma}{2^\sigma-1}$ with values in $G^\delta_{\sigma,s}(\mathbb{R})$. Moreover $T_0\approx\frac{1}{\|u_0\|_{G^1_{\sigma,s}(\mathbb{R})}}$.
\begin{proof}
In order to use Theorem \ref{Th1}, we rewrite (CH) as follows:
\begin{equation}
\left\{
\begin{array}{ll}
u_t=F(u)\doteq-uP_3u-P_{13}[u^2+\frac{1}{2}(P_3u)^2],\\[1ex]
u|_{t=0}=u_0.
\end{array}
\right.
\end{equation}
For a fixed $\sigma\geq1$ and $s>\frac{3}{2}$. By virtue of Proposition \ref{p3}, we have $\{G^\delta_{\sigma,s}\}_{0<\delta<1}$ is a scale of decreasing Banach spaces.
Let $C_s$ be the constant given in Proposition \ref{Product}. By virtue of  Propositions \ref{D}, \ref{Product} and \ref{Mult}, we deduce that for any $0<\delta'<\delta$,
\begin{align}
\|F(u)\|_{G^{\delta'}_{\sigma,s}}&\leq \frac{1}{2}\|P_3(u^2)\|_{G^{\delta'}_{\sigma,s}}+\frac{1}{2}\|u^2\|_{G^{\delta'}_{\sigma,s}}+\frac{1}{2}\|(P_3u)^2\|_{G^{\delta'}_{\sigma,s-1}}\\
\nonumber&\leq C_s\frac{e^{-\sigma}\sigma^{\sigma}}{2(\delta-\delta')^\sigma}\|u\|^2_{G^\delta_{\sigma,s}}+\frac{C_s}{2}\|u\|^2_{G^\delta_{\sigma,s}}+\frac{C_s}{2}\|P_3u\|^2_{G^\delta_{\sigma,s-1}}\leq \frac{C_s (e^{-\sigma}\sigma^{\sigma}+2)}{2(\delta-\delta')^\sigma}\|u\|^2_{G^\delta_{\sigma,s}},
\end{align}
which implies that $F$ satisfies the condition (1) of Theorem \ref{Th1}. By the same token, we obtain that $\|F(u_0)\|_{G^\delta_{\sigma,s}}\leq \frac{C_s (e^{-\sigma}\sigma^{\sigma}+2)}{2(1-\delta)^\sigma}\|u_0\|^2_{G^1_{\sigma,s}}$. Thus, we see that $F$ satisfies the condition (3) of Theorem \ref{Th1} with $M=C_s(\frac{e^{-\sigma}\sigma^{\sigma}}{2}+1)\|u_0\|^2_{G^1_{\sigma,s}}$. In order to prove our desire result, it suffices to show that $F$ satisfies the condition (2) of Theorem \ref{Th1}. Assume that $\|u-u_0\|_{G^\delta_{\sigma,s}}\leq R$ and $\|v-u_0\|_{G^\delta_{\sigma,s}}\leq R$. Applying Propositions \ref{D} and \ref{Mult}, we get
\begin{align}
\|F(u)-F(v)\|_{G^{\delta'}_{\sigma,s}}&\leq \frac{e^{-\sigma}\sigma^{\sigma}}{2(\delta-\delta')^\sigma}\|u^2-v^2\|_{G^{\delta}_{\sigma,s}}+\|P_{13}(u^2-v^2)\|_{G^{\delta'}_{\sigma,s}}+\frac{1}{2}\|P_{13}[(P_3u)^2-(P_3v)^2]\|_{G^{\delta'}_{\sigma,s}}\\
\nonumber&\leq  \frac{e^{-\sigma}\sigma^{\sigma}}{2(\delta-\delta')^\sigma} \|u^2-v^2\|_{G^{\delta}_{\sigma,s}}+\frac{1}{2}\|u^2-v^2\|_{G^{\delta}_{\sigma,s}}+\frac{1}{2} \|(P_3u)^2-(P_3v)^2\|_{G^{\delta}_{\sigma,s-1}}\\
\nonumber&\leq \frac{C_s(e^{-\sigma}\sigma^{\sigma}+1)}{2(\delta-\delta')^\sigma}\|u+v\|_{G^{\delta}_{\sigma,s}}\|u-v\|_{G^{\delta}_{\sigma,s}}+\frac{C_s}{2(\delta-\delta')^\sigma}\|P_3u+P_3v\|_{G^{\delta}_{\sigma,s-1}}\|P_3u-P_3v\|_{G^{\delta}_{\sigma,s-1}}\\
\nonumber&\leq \frac{C_s(e^{-\sigma}\sigma^{\sigma}+1)}{2(\delta-\delta')^\sigma}\|u+v\|_{G^{\delta}_{\sigma,s}}\|u-v\|_{G^{\delta}_{\sigma,s}}+\frac{C_s}{2(\delta-\delta')^\sigma}\|u+v\|_{G^{\delta}_{\sigma,s}}\|u-v\|_{G^{\delta}_{\sigma,s}}\\
\nonumber&\leq \frac{C_s(e^{-\sigma}\sigma^{\sigma}+2)}{2(\delta-\delta')^\sigma}(\|u\|_{G^{\delta}_{\sigma,s}}+\|v\|_{G^{\delta}_{\sigma,s}})\|u-v\|_{G^{\delta}_{\sigma,s}}\\
\nonumber&\leq \frac{C_s(e^{-\sigma}\sigma^{\sigma}+2)}{(\delta-\delta')^\sigma}(\|u_0\|_{G^{\delta}_{\sigma,s}}+R)\|u-v\|_{G^{\delta}_{\sigma,s}}\leq \frac{C_s(e^{-\sigma}\sigma^{\sigma}+2)}{(\delta-\delta')^\sigma}(\|u_0\|_{G^{1}_{\sigma,s}}+R)\|u-v\|_{G^{\delta}_{\sigma,s}}.
\end{align}
From the above inequality, we verify that $F$ satisfies the condition (2) of Theorem \ref{Th1} with $L=C_s(e^{-\sigma}\sigma^{\sigma}+2)(\|u_0\|_{G^{1}_{\sigma,s}}+R)$. Moreover, $T_0=\min\{\frac{1}{2^{2\sigma+4}L}, \frac{(2^\sigma-1)R}{(2^\sigma-1)2^{2\sigma+3}LR+M}\}$, by setting $R=\|u_0\|_{G^{1}_{\sigma,s}}$, we see that $L=2C_s(e^{-\sigma}\sigma^{\sigma}+2)\|u_0\|_{G^{1}_{\sigma,s}}$ and $M\leq 2^{2\sigma+3}LR$. Then, we have $T_0=\frac{1}{2^{2\sigma+5}C_s(e^{-\sigma}\sigma^{\sigma}+2)\|u_0\|_{G^{1}_{\sigma,s}}}$.
\end{proof}
\end{theo}
\begin{theo}\label{2CHTh}
Let $\sigma\geq1$ and $s>\frac{3}{2}$. Assume that $u_0\in G^1_{\sigma,s}(\mathbb{R})$ and $\rho_0\in G^1_{\sigma,s-1}(\mathbb{R})$. Then for every $0<\delta<1$,  there exists a $T_0>0$ such that the two-component Camassa-Holm system has a unique solution $(u,\rho)$ which is holomorphic in $|t|<\frac{T_0(1-\delta)^\sigma}{2^\sigma-1}$ with values in $G^\delta_{\sigma,s}(\mathbb{R})\times G^\delta_{\sigma,s-1}(\mathbb{R})$. Moreover $T_0\approx\frac{1}{\|u_0\|_{G^1_{\sigma,s}(\mathbb{R})}+\|\rho_0\|_{G^1_{\sigma,s-1}(\mathbb{R})}}$.
\begin{proof}
We only consider the case $k=1$, and change the 2-component Camassa-Holm (2CH) system into the following form
  \begin{align}
\left\{
\begin{array}{ll}
z_t=F(z),     \\[1ex]
z|_{t=0}=z_0,
\end{array}
\right.
\end{align}
where $z=(u,\rho)^{\mathrm{T}}$, $z_0=(u_0,\rho_0)^{\mathrm{T}}$ and
\begin{equation}
F(z)=\begin{pmatrix}
F_1(z) \\
F_2(z)
\end{pmatrix}
=\begin{pmatrix}
-P_3(\frac{u^2}{2})-P_{13}[u^2+\frac{1}{2}(P_3u)^2+\frac{1}{2}\rho^2]\\
-P_3(u\rho)
\end{pmatrix}.
\end{equation}
For fixed $\sigma\geq 1$ and $s>\frac{3}{2}$, we set $X_{\delta}=G^{\delta}_{\sigma,s}(\mathbb{R})\times G^{\delta}_{\sigma,s-1}(\mathbb{R})$ and
$$\|z\|_{\delta}=\|u\|_{G^\delta_{\sigma,s}}+\|\rho\|_{G^\delta_{\sigma,s-1}}.$$
Thanks to Proposition \ref{p3}, we have $\{X_\delta\}_{0<\delta<1}$ is a scale of decreasing Banach spaces. Let $C_s$ be the constant given in Proposition \ref{Product}. According to Propositions \ref{D}, \ref{Product} and \ref{Mult}, we have for any $0<\delta'<\delta$,
\begin{align}
\|F_1(z)\|_{G^{\delta'}_{\sigma,s}}&\leq \frac{1}{2}\|P_3(u^2)\|_{G^{\delta'}_{\sigma,s}}+\frac{1}{2}\|u^2\|_{G^{\delta'}_{\sigma,s}}+\frac{1}{2}\|(P_3u)^2\|_{G^{\delta'}_{\sigma,s-1}}+\frac{1}{2}\|\rho^2\|_{G^{\delta'}_{\sigma,s-1}}\\
\nonumber&\leq C_s\frac{e^{-\sigma}\sigma^{\sigma}}{2(\delta-\delta')^\sigma}\|u\|^2_{G^\delta_{\sigma,s}}+\frac{C_s}{2}\|u\|^2_{G^\delta_{\sigma,s}}+\frac{C_s}{2}\|P_3u\|^2_{G^\delta_{\sigma,s-1}}+\frac{C_s}{2}\|\rho\|^2_{G^\delta_{\sigma,s-1}}
\\
\nonumber&\leq \frac{C_s (e^{-\sigma}\sigma^{\sigma}+2)}{2(\delta-\delta')^\sigma}\|u\|^2_{G^\delta_{\sigma,s}}+\frac{C_s}{2(\delta-\delta')^\sigma}\|\rho\|^2_{G^\delta_{\sigma,s-1}}\leq \frac{C_s (e^{-\sigma}\sigma^{\sigma}+3)}{2(\delta-\delta')^\sigma}\|z\|^2_\delta,
\end{align}
\begin{align}
\|F_2(z)\|_{G^{\delta'}_{\sigma,s-1}}\leq \frac{e^{-\sigma}\sigma^{\sigma}}{(\delta-\delta')^\sigma}\|u\rho\|_{G^{\delta}_{\sigma,s-1}}\leq \frac{C_se^{-\sigma}\sigma^{\sigma}}{(\delta-\delta')^\sigma}\|u\|_{G^{\delta}_{\sigma,s}}\|\rho\|_{G^{\delta}_{\sigma,s-1}}\leq \frac{C_se^{-\sigma}\sigma^{\sigma}}{(\delta-\delta')^\sigma}\|z\|^2_{\delta},
\end{align}
which imply that $\|F(z)\|_{\delta'}=\|F_1(z)\|_{G^{\delta'}_{\sigma,s}}+\|F_2(z)\|_{G^{\delta'}_{\sigma,s-1}}\leq \frac{C_s (e^{-\sigma}\sigma^{\sigma}+5)}{2(\delta-\delta')^\sigma}\|z\|^2_\delta$ and $F$ satisfies the condition (1) of Theorem \ref{Th1}. By the same token, we obtain that $\|F(z_0)\|_{\delta}\leq \frac{C_s (e^{-\sigma}\sigma^{\sigma}+5)}{2(1-\delta)^\sigma}\|u_0\|^2_{G^1_{\sigma,s}}$. Thus, we see that $F$ satisfies the condition (3) of Theorem \ref{Th1} with $M=\frac{C_s(e^{-\sigma}\sigma^{\sigma}+5)}{2}\|z_0\|^2_{1}$. In order to prove our desire result, it suffices to show that $F$ satisfies the condition (2) of Theorem \ref{Th1}. Assume that $\|z_1-z_0\|_{\delta}\leq R$ and $\|z_2-z_0\|_{\delta}\leq R$. Taking advantage of Propositions \ref{D} and \ref{Mult}, we get
\begin{align}
\|F_1(z_1)-F_1(z_2)\|_{G^{\delta'}_{\sigma,s}}&\leq \frac{e^{-\sigma}\sigma^{\sigma}}{2(\delta-\delta')^\sigma}\|u^2_1-u^2_2\|_{G^{\delta}_{\sigma,s}}+\|P_{13}(u^2_1-u^2_2)\|_{G^{\delta'}_{\sigma,s}}
\\
\nonumber&+\frac{1}{2}\|P_{13}[(P_3u_1)^2-(P_3u_2)^2]\|_{G^{\delta'}_{\sigma,s}}
+\frac{1}{2}\|P_{13}(\rho^2_1-\rho^2_2)\|_{G^{\delta'}_{\sigma,s}}\\
\nonumber&\leq  \frac{e^{-\sigma}\sigma^{\sigma}}{2(\delta-\delta')^\sigma} \|u^2_1-u^2_2\|_{G^{\delta}_{\sigma,s}}+\frac{1}{2}\|u^2_1-u^2_2\|_{G^{\delta}_{\sigma,s}}
\\
\nonumber&+\frac{1}{2} \|(P_3u_1)^2-(P_3u_2)^2\|_{G^{\delta}_{\sigma,s-1}}+\frac{1}{2}\|\rho^2_1-\rho^2_2\|_{G^{\delta}_{\sigma,s-1}}\\
\nonumber&\leq \frac{C_s(e^{-\sigma}\sigma^{\sigma}+1)}{2(\delta-\delta')^\sigma}\|u_1+u_2\|_{G^{\delta}_{\sigma,s}}\|u_1-u_2\|_{G^{\delta}_{\sigma,s}}\\
\nonumber&+\frac{C_s}{2(\delta-\delta')^\sigma}\|P_3u_1+P_3u_2\|_{G^{\delta}_{\sigma,s-1}}\|P_3u_1-P_3u_2\|_{G^{\delta}_{\sigma,s-1}}
+\frac{C_s}{2}\|\rho_1+\rho_2\|_{G^{\delta}_{\sigma,s-1}}\|\rho_1-\rho_2\|_{G^{\delta}_{\sigma,s-1}}\\
\nonumber&\leq \frac{C_s(e^{-\sigma}\sigma^{\sigma}+1)}{2(\delta-\delta')^\sigma}\|u_1+u_2\|_{G^{\delta}_{\sigma,s}}\|u_1-u_2\|_{G^{\delta}_{\sigma,s}}+\frac{C_s}{2(\delta-\delta')^\sigma}\|u_1+u_2\|_{G^{\delta}_{\sigma,s}}\|u_1-u_2\|_{G^{\delta}_{\sigma,s}}\\
\nonumber&+\frac{C_s}{2}\|\rho_1+\rho_2\|_{G^{\delta}_{\sigma,s-1}}\|\rho_1-\rho_2\|_{G^{\delta}_{\sigma,s-1}}\\
\nonumber&\leq \frac{C_s(e^{-\sigma}\sigma^{\sigma}+2)}{2(\delta-\delta')^\sigma}(\|z_1\|_{\delta}+\|z_2\|_{\delta})\|z_1-z_2\|_{\delta}\\
\nonumber&\leq \frac{C_s(e^{-\sigma}\sigma^{\sigma}+2)}{(\delta-\delta')^\sigma}(\|z_0\|_{\delta}+R)\|z_1-z_2\|_{\delta}\leq \frac{C_s(e^{-\sigma}\sigma^{\sigma}+2)}{(\delta-\delta')^\sigma}(\|z_0\|_{1}+R)\|z_1-z_2\|_{\delta},
\end{align}
\begin{align}
\|F_2(z_1)-F_2(z_2)\|_{G^{\delta'}_{\sigma,s-1}}&\leq \frac{e^{-\sigma}\sigma^{\sigma}}{(\delta-\delta')^\sigma}\|u_1\rho_1-u_2\rho_2\|_{G^{\delta}_{\sigma,s-1}}\leq \frac{e^{-\sigma}\sigma^{\sigma}}{(\delta-\delta')^\sigma}[\|(u_1-u_2)\rho_1\|_{G^{\delta}_{\sigma,s-1}}+\|(\rho_1-\rho_2)u_2\|_{G^{\delta}_{\sigma,s-1}}]\\ \nonumber&\leq\frac{C_se^{-\sigma}\sigma^{\sigma}}{(\delta-\delta')^\sigma}(\|u_1-u_2\|_{G^{\delta}_{\sigma,s}}\|\rho_1\|_{G^{\delta}_{\sigma,s-1}}+\|u_2\|_{G^{\delta}_{\sigma,s}}\|\rho_1-\rho_2\|_{G^{\delta}_{\sigma,s-1}})\\
\nonumber&\leq \frac{C_se^{-\sigma}\sigma^{\sigma}}{(\delta-\delta')^\sigma}(\|z_1\|_{\delta}\|u_1-u_2\|_{G^{\delta}_{\sigma,s}}+\|z_2\|_{\delta}\|\rho_1-\rho_2\|_{G^{\delta}_{\sigma,s-1}})\\
\nonumber&\leq \frac{C_se^{-\sigma}\sigma^{\sigma}}{(\delta-\delta')^\sigma}(\|z_0\|_\delta+R)(\|u_1-u_2\|_{G^{\delta}_{\sigma,s}}+\|\rho_1-\rho_2\|_{G^{\delta}_{\sigma,s-1}})\\
\nonumber&\leq \frac{C_se^{-\sigma}\sigma^{\sigma}}{(\delta-\delta')^\sigma}(\|z_0\|_1+R)\|z_1-z_2\|_{\delta}.
\end{align}
From the above inequalities, we verify that $\|F(z_1)-F(z_2)\|_{\delta'}=\|F_1(z_1)-F_1(z_2)\|_{G^{\delta'}_{\sigma,s}}+\|F_2(z_1)-F_2(z_2)\|_{G^{\delta'}_{\sigma,s-1}}\leq\frac{2C_s(e^{-\sigma}\sigma^{\sigma}+1)}{(\delta-\delta')^\sigma}(\|z_0\|_{1}+R)\|z_1-z_2\|_{\delta}$ and $F$ satisfies the condition (2) of Theorem \ref{Th1} with $L=2C_s(e^{-\sigma}\sigma^{\sigma}+1)(\|z_0\|_1+R)$. Moreover, $T_0=\min\{\frac{1}{2^{2\sigma+4}L}, \frac{(2^\sigma-1)R}{(2^\sigma-1)2^{2\sigma+3}LR+M}\}$, by setting $R=\|z_0\|_1$, we see that $L=4C_s(e^{-\sigma}\sigma^{\sigma}+1)\|z_0\|_1$ and $M\leq 2^{2\sigma+3}LR$. Then, we get that $T_0=\frac{1}{2^{2\sigma+6}C_s(e^{-\sigma}\sigma^{\sigma}+1)\|z_0\|_1}$.
\end{proof}
\end{theo}
\begin{rema}
By the similar argument as in the proof of the above theorem, one can obtain the Gevrey regularity and analyticity for the modify 2-component Camassa-Holm system (M2CH).
\end{rema}
\begin{theo}\label{MTh}
Let $\sigma\geq1$ and $s>\frac{1}{2}$. Assume that $(u_0,v_0,w_0)\in (G^1_{\sigma,s}(\mathbb{R}))^3$. Then for every $0<\delta<1$,  there exists a $T_0>0$ such that the three-component Camassa-Holm system has a unique solution $(u,v,w)$ which is holomorphic in $|t|<\frac{T_0(1-\delta)^\sigma}{2^\sigma-1}$ with values in $(G^\delta_{\sigma,s}(\mathbb{R}))^3$. Moreover $T_0\approx\frac{1}{(\|u_0\|_{G^1_{\sigma,s}}+\|v_0\|_{G^1_{\sigma,s}}+\|w_0\|_{G^1_{\sigma,s}})^2+\|u_0\|_{G^1_{\sigma,s}}+\|v_0\|_{G^1_{\sigma,s}}+\|w_0\|_{G^1_{\sigma,s}}}$.
\begin{proof}
By virtue of (3CH), we see that $a=P_1u, c=P_1w$ and $b=P_2(w\cdot P_{13}u-u\cdot P_{13}w)+2P_2(P_{13}u\cdot P_1w-P_1u\cdot P_{13}w)-2P_2v=B(u,w)-2P_2v$. Hence, we change (3CH) into
  \begin{align}
\left\{
\begin{array}{ll}
U_t=F(U),     \\[1ex]
U|_{t=0}=U_0,
\end{array}
\right.
\end{align}
where $U=(u,v,w)^{\mathrm{T}}$, $U_0=(u_0,v_0,w_0)^{\mathrm{T}}$ and
\begin{multline}
F(U)=\begin{pmatrix}
F_1(U) \\
F_2(U) \\
F_3(U)
\end{pmatrix}\\
=\begin{pmatrix}
-v\cdot P_{13}u+P_3u( B(u,w)-2 P_2v)+\frac{3}{2}u (P_3B(u,w)-2P_{23}v)-\frac{3}{2}u(P_{13}u\cdot P_{13}w-P_1u\cdot P_1w)\\
2v\cdot P_3B(u,w)-4vP_{23}v+P_3v\cdot B(u,w)-2P_3vP_2v \\
-v\cdot P_{13}w+P_3w (B(u,w)-2P_2v)+\frac{3}{2}w (P_3B(u,w)- 2P_{23}v)+\frac{3}{2}w(P_{13}u\cdot P_{13}w-P_1u\cdot P_1w)
\end{pmatrix}.
\end{multline}
For fixed $\sigma\geq 1$ and $s>\frac{1}{2}$, we set $X_{\delta}=(G^{\delta}_{\sigma,s}(\mathbb{R}))^3$ and
$$\|U\|_{\delta}=\|u\|_{G^\delta_{\sigma,s}}+\|v\|_{G^\delta_{\sigma,s}}+\|w\|_{G^\delta_{\sigma,s}}.$$
Due to Proposition \ref{p3}, we have $\{X_\delta\}_{0<\delta<1}$ is a scale of decreasing Banach spaces. \\
Let $C_s$ be the constant given in Proposition \ref{Product}. Taking advantage of Propositions \ref{D}, \ref{Product} and \ref{Mult}, we verify that for any $0<\delta'<\delta$
\begin{multline}
\|F_1(U)\|_{G^{\delta'}_{\sigma,s}}\leq \frac{C_s}{2}\|v\|_{G^{\delta'}_{\sigma,s}}\|u\|_{G^{\delta'}_{\sigma,s}}+\frac{C_se^{-\sigma}\sigma^{\sigma}}{(\delta-\delta')^\sigma}\|u\|_{G^{\delta}_{\sigma,s}}(\|B(u,w)\|_{G^{\delta'}_{\sigma,s}}+\frac{1}{2}\|v\|_{G^{\delta'}_{\sigma,s}})\\
+\frac{3}{2}C_s\|u\|_{G^{\delta'}_{\sigma,s}}(\|P_3B(u,w)\|_{G^{\delta'}_{\sigma,s}}+\frac{1}{2}\|v\|_{G^{\delta'}_{\sigma,s}})+\frac{15}{8}C^2_s\|u\|^2_{G^{\delta'}_{\sigma,s}}\|w\|_{G^{\delta'}_{\sigma,s}}.
\end{multline}
Since $B(u,w)=P_2(w\cdot P_{13}u-u\cdot P_{13}w)+2P_2(P_{13}u\cdot P_1w-P_1u\cdot P_{13}w)$, it follows that
\begin{align}
\|B(u,w)\|_{G^{\delta'}_{\sigma,s}}, ~\|P_3B(u,w)\|_{G^{\delta'}_{\sigma,s}} \leq \frac{3}{4}C_s\|u\|_{G^{\delta'}_{\sigma,s}}\|w\|_{G^{\delta'}_{\sigma,s}}.
\end{align}
Plugging (4.7) into (4.6) yields that
\begin{align}
\|F_1(U)\|_{G^{\delta'}_{\sigma,s}}&\leq \frac{5}{4}C_s\|u\|_{G^{\delta'}_{\sigma,s}}\|v\|_{G^{\delta'}_{\sigma,s}}+3C^2_s\|u\|^2_{G^{\delta'}_{\sigma,s}}\|w\|_{G^{\delta'}_{\sigma,s}}+\frac{C_se^{-\sigma}\sigma^\sigma}{(\delta-\delta')^\sigma}\|u\|_{G^{\delta}_{\sigma,s}}(\frac{3}{4}C_s\|u\|_{G^{\delta'}_{\sigma,s}}\|w\|_{G^{\delta'}_{\sigma,s}}+\frac{1}{2}\|v\|_{G^{\delta'}_{\sigma,s}})\\
\nonumber&\leq \frac{C^2_s}{(\delta-\delta')^\sigma}(3+\frac{3 e^{-\sigma}\sigma^\sigma}{4})\|u\|^2_{G^{\delta}_{\sigma,s}}\|w\|_{G^{\delta}_{\sigma,s}}+\frac{C_s}{(\delta-\delta')^\sigma}(\frac{7}{4}+\frac{e^{-\sigma}\sigma^\sigma}{2})\|u\|_{G^{\delta}_{\sigma,s}}\|v\|_{G^{\delta}_{\sigma,s}})\\
\nonumber&\leq \frac{C_s\|U\|^2_{\delta}}{(\delta-\delta')^\sigma}\bigg[C_s\|U\|_{\delta}(3+\frac{3 e^{-\sigma}\sigma^\sigma}{4})+(\frac{5}{4}+\frac{e^{-\sigma}\sigma^\sigma}{2})\bigg].
\end{align}
By the same token, we have
\begin{align}
\|F_2(U)\|_{G^{\delta'}_{\sigma,s}}\leq \frac{C_s\|U\|^2_{\delta}}{(\delta-\delta')^\sigma}\bigg[C_s\|U\|_{\delta}(\frac{3}{2}+\frac{3 e^{-\sigma}\sigma^\sigma}{4})+(1+\frac{e^{-\sigma}\sigma^\sigma}{2})\bigg],\\
\|F_3(U)\|_{G^{\delta'}_{\sigma,s}}\leq \frac{C_s\|U\|^2_{\delta}}{(\delta-\delta')^\sigma}\bigg[C_s\|U\|_{\delta}(3+\frac{3 e^{-\sigma}\sigma^\sigma}{4})+(\frac{5}{4}+\frac{e^{-\sigma}\sigma^\sigma}{2})\bigg],
\end{align}
which lead to
\begin{align}
\|F(U)\|_{\delta'}=\|F_1(U)\|_{G^{\delta'}_{\sigma,s}}+\|F_2(U)\|_{G^{\delta'}_{\sigma,s}}+\|F_3(U)\|_{G^{\delta'}_{\sigma,s}}\leq  \frac{C_s\|U\|^2_{\delta}}{(\delta-\delta')^\sigma}\bigg[C_s\|U\|_{\delta}(\frac{15}{2}+\frac{9 e^{-\sigma}\sigma^\sigma}{4})+(\frac{9}{2}+\frac{3e^{-\sigma}\sigma^\sigma}{2})\bigg].
\end{align}
Thus, we verify that $F$ satisfies the condition (1) of Theorem \ref{Th1}. By the similar estimates as above, we obtain that
\begin{align}
\|F(U_0)\|_{\delta}\leq \frac{C_s\|U_0\|^2_{1}}{(1-\delta)^\sigma}\bigg[C_s\|U_0\|_{1}(\frac{15}{2}+\frac{9 e^{-\sigma}\sigma^\sigma}{4})+(\frac{9}{2}+\frac{3e^{-\sigma}\sigma^\sigma}{2})\bigg],
\end{align}
which implies that $F$ satisfies the condition (3) of Theorem \ref{Th1} and $M=C_s\|U_0\|^2_{1}\bigg[C_s\|U_0\|_{1}(\frac{15}{2}+\frac{9 e^{-\sigma}\sigma^\sigma}{4})+(\frac{9}{2}+\frac{3e^{-\sigma}\sigma^\sigma}{2})\bigg].$ In order to prove our desire result, it suffices to show that $F$ satisfies the condition (2) of Theorem \ref{Th1}. Assume that $\|U_1-U_0\|_{\delta}\leq R$ and $\|U_2-U_0\|_{\delta}\leq R$ for every $0<\delta<1$. Making use of Propositions \ref{Product} and \ref{Mult}, we obtain that for any $0<\delta'<\delta<1$
\begin{align}
\|I\|_{G^{\delta'}_{\sigma,s}}\doteq\|v_1P_{13}u_1-v_2P_{13}u_2\|_{G^{\delta'}_{\sigma,s}}&\leq \frac{1}{2}C_s\|v_1-v_2\|_{G^{\delta'}_{\sigma,s}}\|u_1\|_{G^{\delta'}_{\sigma,s}}+\frac{1}{2}C_s\|u_1-u_2\|_{G^{\delta'}_{\sigma,s}}\|v_2\|_{G^{\delta'}_{\sigma,s}}\\
 \nonumber&\leq C_s\|U_1-U_2\|_{\delta}(\|U_0\|_{1}+R)\leq\frac{C_s(\|U_0\|_{1}+R)}{(\delta-\delta')^\sigma}\|U_1-U_2\|_{\delta},
\end{align}
\begin{align}
\|II\|_{G^{\delta'}_{\sigma,s}}\doteq\|P_3u_1P_2v_1-P_3u_2P_2v_2\|_{G^{\delta'}_{\sigma,s}}&\leq \frac{C_se^{-\sigma}\sigma^{\sigma}}{4(\delta-\delta')^\sigma}\|u_1-u_2\|_{G^\delta_{\sigma,s}}\|v_1\|_{G^\delta_{\sigma,s}}+\frac{C_se^{-\sigma}\sigma^{\sigma}}{4(\delta-\delta')^\sigma}\|u_2\|_{G^\delta_{\sigma,s}}\|v_1-v_2\|_{G^\delta_{\sigma,s}}\\
\nonumber&\leq \frac{C_s(\|U_0\|_1+R)e^{-\sigma}\sigma^{\sigma}}{2(\delta-\delta')^\sigma}\|U_1-U_2\|_{\delta},
\end{align}
\begin{align}
\|III\|_{G^{\delta'}_{\sigma,s}}\doteq\|u_1P_{23}v_1-u_2P_{23}v_2\|_{G^{\delta'}_{\sigma,s}}&\leq  \frac{1}{4}C_s\|v_1-v_2\|_{G^{\delta'}_{\sigma,s}}\|u_1\|_{G^{\delta'}_{\sigma,s}}+\frac{1}{4}C_s\|u_1-u_2\|_{G^{\delta'}_{\sigma,s}}\|v_2\|_{G^{\delta'}_{\sigma,s}}\\
 \nonumber&\leq \frac{C_s}{2}\|U_1-U_2\|_{\delta}(\|U_0\|_{1}+R)\leq\frac{C_s(\|U_0\|_{1}+R)}{2(\delta-\delta')^\sigma}\|U_1-U_2\|_{\delta},
\end{align}
\begin{align}
\|IV\|_{G^{\delta'}_{\sigma,s}}&\doteq\|u_1P_{13}u_1P_{13}w_1-u_2P_{13}u_2P_{13}w_2\|_{G^{\delta'}_{\sigma,s}}\leq\frac{C^2_s}{4}\|u_1-u_2\|_{G^{\delta'}_{\sigma,s}}\|u_1\|_{G^{\delta'}_{\sigma,s}}\|w_1\|_{G^{\delta'}_{\sigma,s}}\\
\nonumber&+\frac{C^2_s}{4}\|u_1-u_2\|_{G^{\delta'}_{\sigma,s}}\|u_2\|_{G^{\delta'}_{\sigma,s}}\|w_1\|_{G^{\delta'}_{\sigma,s}}
+\frac{C^2_s}{4}\|w_1-w_2\|_{G^{\delta'}_{\sigma,s}}\|u_2\|^2_{G^{\delta'}_{\sigma,s}}
\leq \frac{3C^2_s(\|U_0\|_1+R)^2}{4(\delta-\delta')^\sigma}\|U_1-U_2\|_{\delta},
\end{align}
\begin{align*}
\|V\|_{G^{\delta'}_{\sigma,s}}\doteq\|u_1P_{1}u_1P_{1}w_1-u_2P_{1}u_2P_{1}w_2\|_{G^{\delta'}_{\sigma,s}}\leq \frac{3C^2_s(\|U_0\|_1+R)^2}{(\delta-\delta')^\sigma}\|U_1-U_2\|_{\delta},
\end{align*}
\begin{align}
\|B(u_1,w_1)-B(u_2,w_2)\|_{G^{\delta'}_{\sigma,s}}&\leq C(P_2)\|w_1P_{13}u_1-w_2P_{13}u_2\|_{G^{\delta'}_{\sigma,s}}+C(P_2) \|u_1P_{13}w_1-u_2P_{13}w_2\|_{G^{\delta'}_{\sigma,s}}\\
\nonumber&+2C(P_2)\|P_{13}u_1P_{1}w_1-P_{13}u_2P_1w_2\|_{G^{\delta'}_{\sigma,s}}+2C(P_2)\|P_{13}w_1P_{1}u_1-P_{13}w_2P_1u_2\|_{G^{\delta'}_{\sigma,s}}\\
\nonumber&\leq \frac{3C_s(\|U_0\|_1+R)}{2(\delta-\delta')^\sigma}\|U_1-U_2\|_{\delta},
\end{align}
\begin{align}
\|P_3B(u_1,w_1)-P_3B(u_2,w_2)\|&\leq C(P_{23})\|w_1P_{13}u_1-w_2P_{13}u_2\|_{G^{\delta'}_{\sigma,s}}+C(P_{23}) \|u_1P_{13}w_1-u_2P_{13}w_2\|_{G^{\delta'}_{\sigma,s}}\\
\nonumber&+2C(P_{23})\|P_{13}u_1P_{1}w_1-P_{13}u_2P_1w_2\|_{G^{\delta'}_{\sigma,s}}+2C(P_{23})\|P_{13}w_1P_{1}u_1-P_{13}w_2P_1u_2\|_{G^{\delta'}_{\sigma,s}}\\
\nonumber&\leq \frac{3C_s(\|U_0\|_1+R)}{2(\delta-\delta')^\sigma}\|U_1-U_2\|_{\delta},
\end{align}
\begin{align}
\|VI\|_{G^{\delta'}_{\sigma,s}}&\doteq\|P_3u_1B(u_1,w_1)-P_3u_2B(u_2,w_2)\|_{G^{\delta'}_{\sigma,s}}\leq \frac{C_se^{-\sigma}\sigma^\sigma}{(\delta-\delta')^\sigma}\|u_1-u_2\|_{G^{\delta}_{\sigma,s}}\|B(u_1,w_1)\|_{G^{\delta'}_{\sigma,s}}\\
\nonumber&+\frac{C_se^{-\sigma}\sigma^\sigma}{(\delta-\delta')^\sigma}\|u_2\|_{G^{\delta}_{\sigma,s}}\|B(u_1,w_1)-B(u_2,w_2)\|_{G^{\delta'}_{\sigma,s}}\leq \frac{9C^2_se^{-\sigma}\sigma^\sigma(\|U_0\|_1+R)^2}{4(\delta-\delta')^\sigma} \|U_1-U_2\|_{\delta},
\end{align}
\begin{align}
\|VII\|_{G^{\delta'}_{\sigma,s}}&\doteq\|u_1P_3B(u_1,w_1)-u_2P_3B(u_2,w_2)\|_{G^{\delta'}_{\sigma,s}}\leq C_s\|u_1-u_2\|_{G^{\delta'}_{\sigma,s}}\|P_3B(u_1,w_1)\|_{G^{\delta'}_{\sigma,s}}
\\
\nonumber&+C_s\|u_2\|_{G^{\delta'}_{\sigma,s}}\|P_3B(u_1,w_1)-P_3B(u_2,w_2)\|_{G^{\delta'}_{\sigma,s}}\leq \frac{9C^2_s(\|U_0\|_1+R)^2}{4(\delta-\delta')^\sigma} \|U_1-U_2\|_{\delta}.
\end{align}
Since $F_1(U_1)-F_1(U_2)=-I-2II-3III-\frac{3}{2}IV+\frac{3}{2}V+VI+\frac{3}{2}VII$, it follows from the above inequalities that
\begin{align}
\|F_1(U_1)-F_1(U_2)\|_{G^{\delta'}_{\sigma,s}}\leq  \frac{C_s\|U_1-U_2\|_{\delta}(\|U_0\|_1+R)}{(\delta-\delta')^\sigma}\bigg[C_s(\|U_0\|_{1}+R)(9+\frac{9 e^{-\sigma}\sigma^\sigma}{4})+(\frac{5}{2}+e^{-\sigma}\sigma^\sigma)\bigg].
\end{align}
By the similar way, we deduce that
\begin{align}
\|F_2(U_1)-F_2(U_2)\|_{G^{\delta'}_{\sigma,s}}\leq  \frac{C_s\|U_1-U_2\|_{\delta}(\|U_0\|_1+R)}{(\delta-\delta')^\sigma}\bigg[C_s(\|U_0\|_{1}+R)(\frac{9}{2}+\frac{9 e^{-\sigma}\sigma^\sigma}{4})+(2+e^{-\sigma}\sigma^\sigma)\bigg],
\end{align}
\begin{align}
\|F_3(U_1)-F_3(U_2)\|_{G^{\delta'}_{\sigma,s}}\leq  \frac{C_s\|U_1-U_2\|_{\delta}(\|U_0\|_1+R)}{(\delta-\delta')^\sigma}\bigg[C_s(\|U_0\|_{1}+R)(9+\frac{9 e^{-\sigma}\sigma^\sigma}{4})+(\frac{5}{2}+e^{-\sigma}\sigma^\sigma)\bigg],
\end{align}
which lead to
\begin{align}
\|F(U_1)-F(U_2)\|_{\delta'}\leq \frac{C_s\|U_1-U_2\|_{\delta}(\|U_0\|_1+R)}{(\delta-\delta')^\sigma}\bigg[C_s(\|U_0\|_{1}+R)(\frac{45}{2}+\frac{27 e^{-\sigma}\sigma^\sigma}{4})+(7+3e^{-\sigma}\sigma^\sigma)\bigg].
\end{align}
The above inequality implies that $F$ satisfies the condition (2) of Theorem \ref{Th1} and $L= C_s(\|U_0\|_1+R)\bigg[C_s(\|U_0\|_{1}+R)(\frac{45}{2}+\frac{27 e^{-\sigma}\sigma^\sigma}{4})+(7+3e^{-\sigma}\sigma^\sigma)\bigg]$. Moreover $T_0=\min\{\frac{1}{2^{2\sigma+4}L}, \frac{(2^\sigma-1)R}{(2^\sigma-1)2^{2\sigma+3}LR+M}\}$, by setting $R=\|U_0\|_{1}$, we see that $L=C^2_s\|U_0\|^2_1C_{1,\sigma}+C_s\|U_0\|_1C_{2,\sigma}$, $M\leq 2^{2\sigma+3}LR$. Then, we deduce that $T_0=\frac{1}{2^{2\sigma+4}L}=\frac{1}{2^{2\sigma+4}[C^2_s\|U_0\|^2_1C_{1,\sigma}+C_s\|U_0\|_1C_{2,\sigma}]}$ where $C_{1,\sigma}=90+27e^{-\sigma}\sigma^\sigma$ and $C_{2,\sigma}=14+6e^{-\sigma}\sigma^\sigma$.
\end{proof}
\end{theo}

\section{Continuity of the data-to-solution map}
In this section, we investigate the continuity of the data-to-solution map for initial data and solutions in Theorems \ref{CHTH}, \ref{2CHTh} and \ref{MTh}. We only prove this for the 3-component Camassa-Holm system (3CH) since it is much complex and the proofs are similar for the other systems.
\begin{theo}\label{STH}
Let $\sigma\geq1$ and $s>\frac{1}{2}$. Assume that $(u_0,v_0,w_0)\in (G^1_{\sigma,s}(\mathbb{R}))^3$. Then the data-to-solution map $(u_0,v_0,w_0) \mapsto(u,v,w)$ of the 3-component Camassa-Holm system is continuous from $(G^1_{\sigma,s}(\mathbb{R}))^3$ into the solutions space.
\end{theo}
Firstly we introduce a definition to explain what means the data-to-solution map is continuous from $(G^1_{\sigma,s}(\mathbb{R}))^3$ into the solutions space.
\begin{defi}
Let $\sigma\geq1$ and $s>\frac{1}{2}$. We say that the data-to-solution map $(u_0,v_0,w_0) \mapsto(u,v,w)$ of the 3-component Camassa-Holm system is continuous if for a given $(u^\infty_0,v^\infty_0,w^\infty_0)\in (G^1_{\sigma,s}(\mathbb{R}))^3$ there exists a  $T=T(\|u^\infty_0\|_{G^1_{\sigma,s}},\|v^\infty_0\|_{G^1_{\sigma,s}},\|w^\infty_0\|_{G^1_{\sigma,s}})>0$ such that for any sequence $(u^n_0,v^n_0,w^n_0)\in (G^1_{\sigma,s}(\mathbb{R}))^3$ and $\|u^n_0-u^\infty_0\|_{G^1_{\sigma,s}}+\|v^n_0-v^\infty_0\|_{G^1_{\sigma,s}}+\|w^n_0-w^\infty_0\|_{G^1_{\sigma,s}}\xrightarrow{n\rightarrow\infty}0$, the corresponding solutions $(u^n,v^n,w^n)$ of (3CH) satisfy $\|u^n-u^\infty\|_{E_T}+\|v^n-v^\infty\|_{E_T}+\|w^n-w^\infty\|_{E_T}\xrightarrow{n\rightarrow\infty}0$, where
\begin{align}
\|u\|_{E_T}=\sup_{|t|<\frac{T(1-\delta)^\sigma}{2^\sigma-1},0<\delta<1}\bigg(\|u(t)\|_{G^\delta_{\sigma,s}}(1-\delta)^\sigma\sqrt{1-\frac{|t|}{T(1-\delta)^\sigma}}\bigg).
\end{align}
\end{defi}
{\bf Proof of Theorem \ref{STH}}. Without loss of generality, we may assume that $t\geq0$. As in the proof of Theorem \ref{MTh}, we use the same notation $U^n=(u^n,v^n,w^n)^{\mathrm{T}}$, $U^n_0=(u^n_0,v^n_0,w^n_0)^{\mathrm{T}}$ and $\|U^n\|_\delta=\|u^n\|_{G^\delta_{\sigma,s}}+\|v^n\|_{G^\delta_{\sigma,s}}+\|w^n\|_{G^\delta_{\sigma,s}}$. Define that
\begin{align}
T^\infty=\frac{1}{2^{2\sigma+4}[C^2_s\|U^\infty_0\|^2_1C_{1,\sigma}+C_s\|U^\infty_0\|_1C_{2,\sigma}]},\quad T^n=\frac{1}{2^{2\sigma+4}[C^2_s\|U^n_0\|^2_1C_{1,\sigma}+C_s\|U^n_0\|_1C_{2,\sigma}]},
\end{align}
where $C_{1,\sigma}=90+27e^{-\sigma}\sigma^\sigma$, $C_{2,\sigma}=14+6e^{-\sigma}\sigma^\sigma$ and $C_s$ is given in Proposition \ref{Product}.  Since $\|U^n_0-U^\infty_0\|_1\xrightarrow{n\rightarrow\infty}0$, it follows that there exists a constant $N$ such that, if $n\geq N$ we have
\begin{align}
\|U^n_0\|_1\leq \|U^\infty_0\|_1+1.
\end{align}
By setting
\begin{align}
T=\frac{1}{2^{2\sigma+4}[C^2_s(\|U^\infty_0\|_1+1)^2C_{1,\sigma}+C_s(\|U^\infty_0\|_1+1)C_{2,\sigma}]},
\end{align}
we deduce from (5.2) that $T<\min\{T^n,T^\infty\}$ for any $n\geq N$. As in the proof of Theorem \ref{MTh}, we see that $T^n$ and $T^\infty$ are the existence time corresponding to $\|U^n_0\|_1$ and $\|U^\infty_0\|$ respectively, which implies that for any $n\geq N$
\begin{align}
&U^\infty(t,x)=U^\infty_0(x)+\int^t_0F(U^\infty(t,x))d\tau, \quad 0\leq t<\frac{T(1-\delta)^\sigma}{(2^\sigma-1)}, \\
&U^n(t,x)=U^n_0(x)+\int^t_0F(U^n(t,x))d\tau, \quad 0\leq t<\frac{T(1-\delta)^\sigma}{(2^\sigma-1)},
\end{align}
where $F$ is given in (4.11). From the above equations, we verify that for any $0\leq t<\frac{T(1-\delta)^\sigma}{(2^\sigma-1)}$ and $0<\delta<1$
\begin{align}
\|U^n(t)-U^\infty(t)\|_\delta\leq \|U^\infty_0-U^n_0\|_\delta+\int^t_0\|F(U^n(\tau))-F(U^\infty(\tau))\|_\delta d\tau.
\end{align}
Define that $\delta(\tau)=\frac{1}{2}(1+\delta)+(\frac{1}{2})^{2+\frac{1}{\sigma}}\bigg\{[(1-\delta)^\sigma-\frac{t}{T}]^\frac{1}{\sigma}-[(1-\delta)^\sigma+(2^{\sigma+1}-1)\frac{t}{T}]^\frac{1}{\sigma}\bigg\} $. By virtue of Lemma \ref{L.2}, we see that $\delta<\delta(\tau)<1$. Taking advantage of (4.30), we obtain $\|F(U^n(\tau))-F(U^\infty(\tau))\|_\delta\leq \frac{L\|U^n(t)-U^\infty(t)\|_{\delta(\tau)}}{(\delta(\tau)-\delta)^\sigma}$ where $L=C^2_s\|U_0\|^2_1C_{1,\sigma}+C_s\|U_0\|_1C_{2,\sigma}$. Plugging it into (5.7) yields that
\begin{align}
\|U^n(t)-U^\infty(t)\|_\delta\leq \|U^\infty_0-U^n_0\|_\delta+L\int^t_0 \frac{\|U^n(t)-U^\infty(t)\|_{\delta(\tau)}}{(\delta(\tau)-\delta)^\sigma}d\tau.
\end{align}
Applying Lemma \ref{L.2} with $a=T$, we deduce that
\begin{align}
\|U^n(t)-U^\infty(t)\|_\delta\leq \|U^\infty_0-U^n_0\|_\delta+L\frac{T2^{2\sigma+3}\|U^n-U^\infty\|_{E_T}}{(1-\delta)^\sigma}\sqrt{\frac{T(1-\delta)^\sigma}{T(1-\delta)^\sigma-t}}.
\end{align}
Since $T=\frac{1}{2^{2\sigma+4}[C^2_s(\|U^\infty_0\|_1+1)^2C_{1,\sigma}+C_s(\|U^\infty_0\|_1+1)C_{2,\sigma}]}$ and $L=C^2_s\|U_0\|^2_1C_{1,\sigma}+C_s\|U_0\|_1C_{2,\sigma}$, it follows that $LT2^{2\sigma+3}<\frac{1}{2}$. Then, we have
\begin{align}
\|U^n(t)-U^\infty(t)\|_\delta\leq \|U^\infty_0-U^n_0\|_\delta+\frac{1}{2(1-\delta)^\sigma}\|U^n-U^\infty\|_{E_T}\sqrt{\frac{T(1-\delta)^\sigma}{T(1-\delta)^\sigma-t}},
\end{align}
which leads to
\begin{align}
\|U^n(t)-U^\infty(t)\|_\delta (1-\delta)^\sigma\sqrt{1-\frac{t}{T(1-\delta)^\sigma}} &\leq  \|U^\infty_0-U^n_0\|_\delta (1-\delta)^\sigma\sqrt{1-\frac{t}{T(1-\delta)^\sigma}} +\frac{1}{2}\|U^n-U^\infty\|_{E_T}\\
\nonumber&\leq \|U^\infty_0-U^n_0\|_1 +\frac{1}{2}\|U^n-U^\infty\|_{E_T}.
\end{align}
Note that the right hand side of the above inequality is independent of $t$ and $\delta$. By taking the supremum over $0<\delta<1, 0<t<\frac{T(1-\delta)^\sigma}{2^\sigma-1}$, we obtain that
\begin{align}
\|U^n-U^\infty\|_{E_T}\leq \|U^\infty_0-U^n_0\|_1 +\frac{1}{2}\|U^n-U^\infty\|_{E_T},
\end{align}
which implies that
 \begin{align}
\|U^n-U^\infty\|_{E_T}\leq 2\|U^\infty_0-U^n_0\|_1.
\end{align}
The above inequality holds true for any $n\geq N$ and leads to our desire result.

\begin{rema}
In the period case, the Sobolev-Gevrey norm can be stated as follows
\begin{align}
\|f\|_{G^\delta_{\sigma,s}(\mathbb{T})}=\bigg(\sum_{k\in\mathbb{Z}}(1+|k|^2)^se^{2\delta|k|^{\frac{1}{\sigma}}}|\widehat{f}(k)|^2\bigg)^{\frac{1}{2}}=\|e^{\delta(-\Delta)^{\frac{1}{2\sigma}}}f\|_{H^s(\mathbb{T})},
\end{align}
and the similar propositions still hold true. Taking advantage of Theorem \ref{Th1} and by virtue of the same argument as in Theorems \ref{CHTH}, \ref{2CHTh}, \ref{MTh} and \ref{STH}, we get the similar Gevrey regularity and analytic for the Camassa-Holm type systems.
\end{rema}
{\bf Acknowledgements}.  This work was
partially supported by NNSFC (No.11271382), RFDP (No.
20120171110014), the Macao Science and Technology Development Fund (No. 098/2013/A3) and the key project of Sun Yat-sen University. 

\end{document}